\documentclass[10pt]{amsart}

\title{Fill Radius and the Fundamental Group}

\author{Mohan Ramachandran}
\address{Department of Mathematics \\ 
                 SUNY at Buffalo \\
                 Buffalo, NY, 14260}

\author{Jon Wolfson}
\address{Department of Mathematics \\ 
                 Michigan State University \\
                 East Lansing, MI 48824}

\thanks{ The second author was 
partially supported by NSF grant DMS-0604759 }

\date{\today}

\newtheorem{thm}{Theorem}[section]
\newtheorem{lem}[thm]{Lemma}

\newtheorem{conj}{Conjecture}[section]
\theoremstyle{definition}

\newtheorem{defn}{Definition}[section]

\numberwithin{equation}{section}

\renewcommand{\a}{\alpha}
\renewcommand{\b}{\beta}
\renewcommand{\d}{\delta}

\newcommand{\g}{\gamma}
\newcommand{\G}{\Gamma}

\newcommand{\p}{\partial}

\newcommand{\s}{\sigma}

\renewcommand{\t}{\tau}

\def\Pb{\ifmmode{\Bbb P}\else{$\Bbb P$}\fi}
\def\Z{\ifmmode{\Bbb Z}\else{$\Bbb Z$}\fi}
\def\Q{\ifmmode{\Bbb Q}\else{$\Bbb Q$}\fi}
\def\C{\ifmmode{\Bbb C}\else{$\Bbb C$}\fi}
\def\R{\ifmmode{\Bbb R}\else{$\Bbb R$}\fi}
\def\H{\ifmmode{\Bbb H}\else{$\Bbb H$}\fi}
\def\S{\ifmmode{S^2}\else{$S^2$}\fi}

\def\fill{\operatorname{ fill rad}}
\def\Hfill{\operatorname{H_1 fill rad}}

\def\Ric{\operatorname{Ric}}

\def\S{\mathcal S}

\usepackage{graphicx} 

\begin{document}

\maketitle

\begin{abstract} 
In this note we relate the geometric notion of fill radius with the fundamental group of the manifold. We prove: Suppose that a closed Riemannian manifold M satisfies the property that its universal cover has bounded fill radius.  Then the fundamental group of M is virtually free. We explain the relevance of this theorem to some conjectures on positive isotropic curvature and 2-positive Ricci curvature.
\end{abstract}

\setcounter{secnumdepth}{1}

\setcounter{section}{0}

\section{\bf Introduction}

Let $(M, g)$ be an $n$-dimensional Riemannian manifold. The notion of {\it fill radius}, introduced in [G1], [G-L], [S-Y],  is a type of ``two-dimensional diameter''. Let $\g$ be a smooth simple closed curve in $M$ which bounds a disk in $M$. Set $N_r(\g) = \{ x \in M: d(x, \g) \leq r \}$. We define the {\it   fill radius of $\g$} to be:
$$
\fill (\g) = \sup \{r: \mbox{dist}(\g, \p M) > r \; \mbox{and} \; \g \mbox{ does not bound a disc in } N_r(\g) \}
$$
We say a Riemannian manifold $(M,g)$ has its {\it fill radius bounded by $C$} if  every smooth simple closed curve $\g$ which bounds a disk in $M$ satisfies,
$$
\fill (\g) \leq C.
$$
Clearly if the diameter of $(M, g)$ is bounded so is its fill radius. In particular if for all $p \in M$, $\Ric(p) \geq \a$, where $\a$ is a positive constant,  then there is a constant $C = C(\a)$ such that the fill radius of $M$ is bounded by $C$. It is an interesting problem to find ``positive curvature conditions'' that imply fill radius bounds.  In [G-L] and [S-Y] versions of the following result on positive scalar curvature and fill radius are proved. (Throughout this introduction, for technical reasons related to the solution of the Plateau problem, if $(M, g)$ is not compact we will assume it is complete, its sectional curvature is bounded above and its injectivity radius is bounded below away from zero, i.e., we will assume that $(M,g)$ has bounded geometry. If $(M,g)$ is a cover of a closed Riemannian manifold then these conditions are satisfied.)

\medskip

\begin{thm}[Gromov-Lawson, Schoen-Yau]
Let $(M,g)$ be a complete Riemannian three manifold  with positive scalar curvature $S$ that satisfies $S \geq \a $, for a constant $\a > 0$. Then if $\g$ is a smooth simple closed curve in $M$ which bounds a disk in $M$:
$$
\fill (\g) \leq  {\sqrt {\frac{8}{3}}} \frac{\pi}{{\sqrt \alpha}}
$$
\end{thm}

\bigskip

We next recall two positive curvature conditions that conjecturally imply fill radius bounds.
We say $(M, g)$ has {\it two-positive Ricci curvature} if at each point $p \in M$ the sum of the two smallest eigenvalues of the Ricci curvature at $p$ is positive. We say that the two-positive Ricci curvature is bounded below by $\alpha$ if the sum of the two smallest eigenvalues is greater than $\alpha$. It has been conjectured by the second author [W] that:

\medskip

\begin{conj}
\label{conj:Ricci}
Let $(M,g)$ be a complete Riemannian $n$-manifold with two-positive Ricci curvature is bounded below by $\alpha$, for a constant $\a>0$. Then if $\g$ is a smooth simple closed curve in $M$ which bounds a disk in $M$:
$$
\fill (\g) \leq C(\a)
$$
\end{conj}

\medskip

We say $(M, g)$ has {\it positive isotropic curvature bounded below by $\a$} if at each point $p \in M$ and for every orthonormal four frame $\{e_1, e_2, e_3, e_4 \}$ the curvature satisfies:
$$
R_{1313} + R_{1414} + R_{2323} + R_{2424} + 2R_{1234} \geq \a.
$$
It has been conjectured (at least, implicitly by Gromov [G2], Fraser [F]) that:

\medskip

\begin{conj}
\label{conj:isotropic}
Let $(M,g)$ be a complete Riemannian $n$-dimensional manifold with positive isotropic curvature bounded below by $\a$, for a constant $\a>0$.  Then if $\g$ is a smooth simple closed curve in $M$ which bounds a disk in $M$:
$$
\fill (\g) \leq C(\a)
$$
\end{conj}

\vspace{1cm}

A group $G$  is said to be  {\it virtually free} if it possesses a finite index subgroup that is a free group. If $G$ is the fundamental group of a manifold $M$ then $G$ is virtually free if some finite cover of $M$ has fundamental group that is a free group.

\bigskip

In this note we prove:

\begin{thm}
\label{maintheorem}
Let $M$ be a closed Riemannian $n$-manifold. Suppose that the universal cover $\pi: {\tilde M} \to M$ is given the Riemannian metric ${\tilde g}$ such that $\pi$ is a local isometry. If $({\tilde M}, {\tilde g})$ has bounded fill radius  then the fundamental group of $M$ is virtually free.
\end{thm}

\bigskip

If Conjecture \ref{conj:Ricci} is true then Theorem \ref{maintheorem} implies that the fundamental group of a closed $n$-manifold with two-positive Ricci curvature is virtually free. If Conjecture \ref{conj:isotropic} is true then Theorem \ref{maintheorem} implies that the fundamental group of a closed $n$-manifold with positive isotropic curvature is virtually free. We remark that, based on the work of Micallef-Wang [M-W], Gromov [G2] and Fraser [F] explicitly conjecture that the fundamental group of a closed $n$-manifold with positive isotropic curvature is virtually free. However, in light of Theorem \ref{maintheorem}, we attribute Conjecture \ref{conj:isotropic} as above.

We are indebted to Bruce Kleiner for pointing out that a homological version of Theorem \ref{maintheorem} can be proved using our techniques. This is outlined in Section 3. The second author wishes to thank Nick Ivanov for useful discussions.

\bigskip

\section{\bf Fill Radius and the Fundamental Group}

In this section we give the proof of Theorem \ref{maintheorem}. Our approach is based on the notion of the number of ends of a group $G$. There are various definitions of this notion. For our purposes the following definition will suffice:

\begin{defn} Given a  group $G$ we define the number of ends, $e(G)$, of $G$ to be the number of topological ends of $\tilde{K}$, where $\tilde{K} \to K$ is a regular covering of the finite simplicial complex $K$ by the simplicial complex $\tilde{K}$ and $G$ is the group of covering transformations.
\end{defn}

In particular, if $G$ is the fundamental group of a closed manifold $N$ then the number of ends of $G$ is the number of ends of the universal cover $\tilde{N}$ of $N$. It is not difficult to show that a group $G$ can have 0,1,2 or infinitely many ends [E].

\medskip

We will need the following three lemmas.

\begin{lem}
\label{lem:generator}
Let $N$ be a closed manifold. Suppose that $N_0 \to N$ is a covering of  $N$ such that $N_0$ has fundamental group $G$ that is finitely generated and has  exactly one end. Let $\g$ be a simple closed curve in $N_0$ that represents an infinite order generator $[\g]$ of $G$. Let $\tilde{N} \to N_0$ be the universal cover and let $\tilde{\g}$ be the lift of $\g$ to $\tilde{N}$. Then the two ends of $\tilde{\g}$ lie in the same end of $\tilde{N}$.
\end{lem}

\begin{proof}
There is a finite simplicial complex $K$ with regular covering  $\tilde{K}$ such that $G$ acts as the group of covering transformations.
There is an imbedding $\imath: K \to N_0$ that induces an epimorphism of fundamental groups. In particular, the generators of $G$ all lie in $K$.  Then there is an imbedding $\tilde{\imath}:  \tilde{K} \to \tilde{N}$. If $B \subset \tilde{N}$ is compact then $\tilde{\imath}^{-1}(B) \subset \tilde{K}$ is compact. 

Let $\g$ be a simple closed curve in $N_0$ that represents an infinite order generator $[\g]$ of $G$. After a homotopy the lift $\tilde{\g}$ can be assumed to lie in $\tilde{K}$. Since $G$ has exactly one end, any two points on $\tilde{\g}$,  not in $\tilde{\imath}^{-1}(B)$, can be joined by a curve $\a$ in $\tilde{K} \setminus \tilde{\imath}^{-1}(B)$. The curve $ \tilde{\imath}(\a)$ then lies in $ \tilde{N} \setminus B$ and joins points on $\tilde{\g}$ not in $B$. Since this is true for any compact set $B$ the conclusion follows. 
\end{proof}

\bigskip

The next lemma is a version of Lemma \ref{lem:generator} for torsion elements  that are sufficiently long.

\begin{lem}
\label{lem:torsion}
Let $N$ be a closed Riemannian manifold. Suppose that $N_0 \to N$ is a covering of  $N$ such that $N_0$ has fundamental group $G$ that is a finitely generated infinite group with  exactly one end. Let $g_i \in G$ be a sequence and suppose that each $g_i$ is represented by a closed curve  $\g_{g_i}$ beginning and ending at $q \in N_0$ such that each curve lies in a fixed compact region $S$. Let $\tilde{N} \to N_0$ be the universal cover and let $\tilde{\g}_{g_i}$ be a lift of $\g_{g_i}$ to $\tilde{N}$. Denote the distance between the endpoints of $\tilde{\g}_{g_i}$ by $d_i$ and suppose that $d_i \to \infty$. Let $x_i$ be a point on $\tilde{\g}_{g_i}$ such that the distance between $x_i$ and each endpoint is at least $\frac{d_i}{2}$. In addition, suppose that all the $x_i$ lie in a small coordinate ball $B_\d(x)$. Then given $B_R(x)$, $R >> \d$, for sufficiently large $i$ the endpoints of  $\tilde{\g}_{g_i}$ can be joined by a path in $\tilde{N} \setminus B_R(x)$. 
\end{lem}

\begin{proof}
For sufficiently large $i$ neither endpoint of $\tilde{\g}_{g_i}$ lies in a relatively compact region of $\tilde{N} \setminus B_R(x)$. Using the same notation as in  Lemma \ref{lem:generator},  after a homotopy the curve $\tilde{\g}_{g_i}$ can be assumed to lie in $\tilde{K}$. Since $G$ has exactly one end the endpoints of $\tilde{\g}_{g_i}$ lie in the same end of  $\tilde{K}$ and therefore can be joined by a  curve $\a$ in $\tilde{K} \setminus \tilde{\imath}^{-1}(B_R(x))$. The curve $ \tilde{\imath}(\a)$ then lies in $ \tilde{N} \setminus B_R(x)$ and joins the endpoints of $\tilde{\g}_{g_i}$.
\end{proof}

\bigskip

\begin{lem}
\label{lem:curves}
Let $M$ be a complete  manifold with finitely generated fundamental group $G$.  Then there is a  compact subset $S$ of $M$ such that every element of $G$ can be represented by a closed curve beginning and ending at $q \in S$ that lies entirely in $S$.
\end{lem}

\begin{proof}
Choose a finite set of generators and represent each generator by a smooth closed curve beginning and ending at $q \in M$. Then each curve lies in a fixed compact set $S$. The result follows. 
\end{proof}

\bigskip

\begin{thm}
\label{thm:one-end}
Let $N$ be a closed Riemannian manifold. Suppose that the universal cover $\tilde{N}$ has the property that the fill radius  of every  simple closed curve is uniformly bounded above. If $G$ is a  finitely generated subgroup of $\pi_1(N)$ then $G$ cannot have exactly one end.
\end{thm}

\begin{proof}
Assume, by way of contradiction, that the subgroup $G$ of $\pi_1(N)$ has exactly one end. Let  $M$ be a covering of $N$ with fundamental group $\pi_1(M)$ isomorphic to $G$. If $G$ contains an element of infinite order the proof is simpler. We begin with this case though, strictly speaking, this is not necessary.

Assume that $G$ contains a generator of infinite order and
denote by $\g$ a minimal geodesic in $M$ that represents this generator. Let $p: \tilde{N} \to M$ be the universal cover and let $\tilde{\g}$ be the geodesic line that is a lift to $\tilde{N}$ of $\g$. Let $x \in \tilde{\g}$ and $B_R(x) \subset \tilde{N}$ be the metric ball of radius $R$, center $x$. Then because $G$ has exactly one end by Lemma \ref{lem:generator} both ends of $\tilde{\g}$ in $\tilde{N} \setminus B_R(x)$ lie in the same end of $\tilde{N}$. The geodesic line $\tilde{\g}$ consist of two geodesic rays $\tilde{\g}_1$ and $\tilde{\g}_2$ beginning at $x$. For $j=1,2$, choose a point $p_j \in \tilde{N} \setminus B_R(x)$ along $\tilde{\g}_j$  and denote the segment of $\tilde{\g}_j$ from $x$ to $p_j$ by $\t_j$. Since $p_1$ and $p_2$ lie in the same end there is a curve $\b \subset \tilde{N} \setminus B_R(x)$ joining $p_1$ and $p_2$. Denote the closed curve $\t_1 \cup \b \cup \t_2$ by $\eta$. Since $\tilde{N}$ is simply connected $\eta$ is null homotopic and  has fill radius greater than $\frac{R}{2}$. For sufficiently large $R$ this contradicts the fill radius bound.

Next assume that $G$ has no elements of infinite order but that $G$ is infinite. Since $G$ is finitely generated there is a point $q \in M$ and a ball $B_r(q) \subset M$ such that every element $g \in G$ can be represented by a closed curve $\g_g$ in $B_r(q)$ beginning and ending at $q$. Let $\tilde{q}$ denote a lift of $q$ and denote by $\tilde{\g}_g$ the lifts of the $\g_g$ to $\tilde{N}$ that begin at $\tilde{q}$. The endpoints of $\tilde{\g}_g$ are the points $\tilde{q}$ and $g(\tilde{q})$. Since $G$ is infinite we can choose a sequence $g_i \in G$ such that the distance $d(\tilde{q}, g_i(\tilde{q}))= d_i \to \infty$. Choose a point $x_i$ on $\tilde{\g}_{g_i}$ such that $d(x_i, \tilde{q}) \geq \frac{d_i}{2}$ and  $d(x_i, g_i(\tilde{q})) \geq \frac{d_i}{2}$. Using the Deck transformations find elements $h_i$ of $G$ that move the points $x_i$ into a fixed fundamental region $U$ containing $\tilde{q}$. Denote the  curves $h_i(\tilde{\g}_{g_i})$ by $\s_i$. Then the endpoints of $\s_i$ remain at least $\frac{d_i}{2}$ distant from $h_i(x_i)$ and are distance $d_i$ from each other. Note that the sequence $\{h_i(x_i)\}$ lies in the compact set $\bar{U} \cap p^{-1}(B_r(q))$. Choosing a subsequence of $\{h_i(x_i)\}$ we can suppose that the sequence $\{h_i(x_i)\}$ converges to $x \in \tilde{N}$ and therefore that $h_i(x_i) \in B_\d(x)$, for some $\d > 0$. Consider the ball $B_R(x) \in \tilde{N}$, where $R >>\d$. Denote the endpoints of $\s_i$ by $y_i$ and $z_i$. For $i$ sufficiently large, $y_i$ and $z_i$ lie outside $B_R(x)$ and do not lie in any relatively compact region of  $\tilde{N} \setminus B_R(x)$. Thus, by Lemma \ref{lem:torsion},  $y_i$ and $z_i$ can be joined by a smooth curve $\a_i$ lying in $\tilde{N} \setminus B_R(x)$.  Join $x$ to $y_i$ by a minimal geodesic $\t_i$ and join $x$ to $z_i$ by a minimal geodesic $\rho_i$. The closed loop $\t_i \cup \a_i \cup \rho_i$ is null homotopic and  has fill radius greater than $\frac{R}{2}$.  For $R$  sufficiently large this contradicts the fill radius bound.

Finally if $G$ is finite then $G$ has zero ends.
\end{proof}

\bigskip

To prove our next result we will use work of Dunwoody [D].  Stallings' Structure theorem [St1] for finitely generated groups with more than one end is formulated in [D] as follows: {\it  Let $G$ be a finitely generated group. Then $e(G) > 1$ if and only if there is a $G$-tree $T$ such that the stabilizer $G_e$ of each edge $e$ is finite and the stabilizer $G_v$ of each vertex $v$ is finitely generated and  $G_v \neq G$}.

\medskip

\begin{defn} A finitely generated group $G$ is said to be {\it accessible} if there is a $G$-tree $T$ such that $G_e$ is finite for each edge $e$ of $T$ and $G_v$ has at most one end for each vertex $v$ of $T$.
\end{defn}

\medskip

Dunwoody's main result in [D] is: {\it A finitely presented group $G$ is accessible.} (also, see [D-D] Chap. 6 Theorem 6.3).

\begin{thm}
\label{thm:virtually-free}
Let $N$ be a closed Riemannian manifold. Suppose that the universal cover $\pi: {\tilde N} \to N$ is given the Riemannian metric ${\tilde g}$ such that $\pi$ is a local isometry. If $({\tilde N}, {\tilde g})$ has fill radius bounded above then the fundamental group $\pi_1(N)$ is virtually free.
\end{thm}

\begin{proof}
By Theorem \ref{thm:one-end}, $G =\pi_1(N)$ has no finitely generated subgroups with exactly one end. Since $\pi_1(N)$ is finitely presented, by Dunwoody's result, it is accessible. Therefore there is a $G$-tree $T$ such that $G_e$ is finite for each edge $e$ of $T$ and $G_v$ is finite for each vertex $v$ of $T$. Then, by [Se] (see Chap. II, Sec. 2.6, Prop. 11), it follows that $G$ is virtually free.
\end{proof}

\medskip

Under a more restrictive condition on the fundamental group a better result is available.

\begin{thm}
\label{thm:virtually-free2}
Let $N$ be a closed Riemannian manifold with torsion-free fundamental group. Suppose that the universal cover $\pi: {\tilde N} \to N$ is given the Riemannian metric ${\tilde g}$ such that $\pi$ is a local isometry. If $({\tilde N}, {\tilde g})$ has fill radius bounded above then the fundamental group $\pi_1(N)$ is free of finite rank.
\end{thm}

\begin{proof}
We use Grushko's Theorem (see [Ma]) and the following theorem of Stallings [St2] (also, [D-D] Chap. 4 Theorem 6.10): {\it  If $G$ is a torsion-free, finitely generated group with infinitely many ends then $G$ is a non-trivial free product}. Applying Stallings' theorem to $G = \pi_1(N)$, we have $G \simeq G_1 \ast G_2$, where each $G_i$ is finitely generated (by Grushko's Theorem) and each $G_i$ has either two or infinitely many ends (by Theorem \ref{thm:one-end}). Then apply Stallings theorem to each $G_i$ with infinitely many ends and iterate. By Grushko's Theorem, this process terminates after finitely many steps resulting in $G \simeq G_1 \ast \dots \ast G_k$, where each $G_i$  is finitely generated and has  two ends. Since a torsion-free, finitely generated group with two ends is infinite cyclic, we conclude that $G= \pi_1(N)$ is a free group of finite rank.
\end{proof}

\bigskip

\section{\bf Homology and Fill Radius}

In this section we describe the analog of the previous results for the notion of {\it homological fill radius}. We will continue to work with $n$-dimensional Riemannian manifolds $(M,g)$ though the results we describe can be formulated for more general spaces. Let $\G$ be a one-cycle which bounds in $M$. Set $N_r(\G) = \{ x \in M: d(x, \G) \leq r \}$. We define the {\it  homological fill radius of $\G$} to be:
$$
\Hfill (\G) = \sup \{r: \mbox{dist}(\G, \p M) > r \; \mbox{and} \; \G \mbox{ does not bound in } N_r(\G) \}
$$
We say a Riemannian manifold $(M,g)$ has its {\it homological fill radius bounded by $C$} if  every one-cycle $\G$ which bounds  in $M$ satisfies,
$$
\Hfill (\G) \leq C.
$$

\bigskip

Our main theorem is:

\begin{thm}
\label{thm:mainhomology}
Let $(X,g)$ be a complete Riemannian $n$-manifold with bounded geometry, with $H_1(X, \Z) = 0$ and that satisfies a  homological fill radius bound. Suppose the group $G$ acts freely, properly discontinuously and co-compactly on $X$. Then $G$ is virtually free. If, in addition, $G$ is torsion free then $G$ is free of finite rank.
\end{thm}

\bigskip

The first step in the proof of the theorem is the analog of Theorem \ref{thm:one-end}.

\begin{thm}
\label{thm:one-end2}
Let $(X,g)$ be a complete Riemannian $n$-manifold with bounded geometry, with $H_1(X, \Z) = 0$ and that satisfies a  homological fill radius bound. Suppose the group $G$ acts freely, properly discontinuously and co-compactly on $X$. If $H$ is a finitely generated, infinite subgroup of $G$ then $H$ cannot have exactly one end.
\end{thm}

\begin{proof}
To begin we observe, without loss of generality, that we can suppose that $G$ acts as a group of isometries. To see this, note that given a metric on $X/G$ it lifts to a complete metric with bounded geometry $\bar{g}$ on $X$. By assumption $(X, g)$ is a complete Riemannian manifold with bounded geometry. Therefore $\bar{g}$ and $g$ are quasi-isometric. The  homological fill radius bound is a quasi-isometry invariant so $(X, \bar{g})$ satisfies this condition.

Suppose, by way of contradiction, that $H$ has exactly one end.
The quotient space $N = X/G$ is a compact Riemannian manifold with $G$ a subgroup of the fundamental group. In particular, $H$ is a finitely generated, infinite subgroup of the fundamental group. As in the proof of Theorem \ref{thm:one-end}, let  $M$ be a covering of $N$ with fundamental group $\pi_1(M)$ isomorphic to $H$. Note that since $H$ is a subgroup of $G$, $X$ is a regular covering space of $M$. Since $H$ is finitely generated there is a point $q \in M$ and a ball $B_r(q) \subset M$ such that every element $h \in H$ can be represented by a closed curve $\g_h$ in $B_r(q)$ beginning and ending at $q$. Let $\tilde{q} \in X$ denote a lift of $q$ and denote by $\tilde{\g}_h$ the lifts of the $\g_h$ to $X$ that begin at $\tilde{q}$. Since $H$ is infinite we can choose a sequence $h_i \in H$ such that the distance $d(\tilde{q}, h_i(\tilde{q}))= d_i \to \infty$. Choose a point $x_i$ on $\tilde{\g}_{g_i}$ such that $d(x_i, \tilde{q}) \geq \frac{d_i}{2}$ and  $d(x_i, h_i(\tilde{q})) \geq \frac{d_i}{2}$.  On $X$ the group $H$ acts as Deck transformations so using the same argument as in the proof of Theorem \ref{thm:one-end} we can suppose that the points $x_i$ lie in a ball $B_\d(x)$ for some $x \in X$. Hence we find a sequence of curves $\s_i$ in $X$ with endpoints $y_i$ and $z_i$ and containing points $x_i$ with the following properties: (i) all the $x_i$  lie in a fixed coordinate ball $B_\d(x)$ for some $x \in X$, (ii) $d(y_i, z_i) = d_i \to \infty$, as $i \to \infty$, (iii) $d(x_i, y_i) \geq \frac{d_i}{2}$ and $d(x_i, z_i) \geq \frac{d_i}{2}$. Given a metric ball $B_R(x) \subset X$ choose $i$ sufficiently large such that  $y_i$ and $z_i$ lie outside $B_R(x)$ and do not lie in any relatively compact region of  $X \setminus B_R(x)$. 

To use the assumption that $H$ has one end we let $K$ be a finite simplicial complex with regular covering  $\tilde{K}$ such that $H$ acts as the group of covering transformations.
There is an imbedding $\imath: K \to M$ that induces an epimorphism of fundamental groups. In particular, the generators of $H$ all lie in $K$. Since $H$ is a subgroup of $G$, there is an imbedding $\tilde{\imath}:  \tilde{K} \to X$. Then, up to homotopy, $\s_i$ is a path in $\tilde{\imath}(\tilde{K})$. For $i$ sufficiently large the endpoints of $\tilde{\imath}^{-1}(\s_i)$ lie outside $\tilde{\imath}^{-1}(B_R(x))$ and do not lie in any relatively compact region of  $\tilde{K} \setminus \tilde{\imath}^{-1}(B_R(x))$. Since $H$ has exactly one end it follows that the endpoints of $\tilde{\imath}^{-1}(\s_i)$ can be joined by a  curve lying in $\tilde{K} \setminus \tilde{\imath}^{-1}(B_R(x))$. Hence the endpoints of $\s_i$,  $y_i$ and $z_i$, can be joined by a smooth curve $\a_i$ lying in in $X \setminus B_R(x)$. Join $x$ to $y_i$ by a minimal geodesic $\t_i$ and join $x$ to $z_i$ by a minimal geodesic $\rho_i$. The closed loop $\G = \t_i \cup \a_i \cup \rho_i$ is a one-cycle in $X$ and since $H_1(X, \Z) = 0$, $\G$ spans a chain. The homological fill radius of $\G$ is greater than $\frac{R}{2}$.  For $R$  sufficiently large this contradicts the homological fill radius bound of $X$.
\end{proof}

\noindent{\it Proof of Theorem \ref{thm:mainhomology}.} The group $G$ need not be finitely presented however Dunwoody's work applies more generally to show that if $G$ is almost finitely presented then $G$  is accessible [D]. A group that acts freely, properly discontinuously and co-compactly on a space $X$ with $H_1(X, \Z_2) = 0$ is almost finitely presented. Therefore $G$ is accessible. By Theorem \ref{thm:one-end2} no finitely generated subgroup of $G$ has exactly one end.  Hence there is a $G$-tree $T$ such that $G_e$ is finite for each edge $e$ of $T$ and $G_v$ is finite for each vertex $v$ of $T$. Then, by [Se], it follows that $G$ is virtually free. 

If $G$ is torsion free the proof is identical to the proof of Theorem \ref{thm:virtually-free2}.


\vspace{2cm}


















\begin{thebibliography}{XXX}

\bibitem[D-D]{dd} Dicks, W. and Dunwoody, M., {\it Groups acting on graphs}, Cambridge University Press, Cambridge, 1989.

\bibitem[D]{d} Dunwoody, M., The accessibility of finitely presented groups, Invent. Math. {\bf 81} (1985) 449-457.


\bibitem[E]{e} Epstein, D., Ends, {\it Topology of 3-manifolds} edited by, M.K. Fort, Jr., Prentice-Hall, 1962, 110-117.

\bibitem[G1]{g1} Gromov, M, Filling Riemannian manifolds, J. Diff. Geom. {\bf 18} (1983) 1- 147.

\bibitem[G2]{g2} Gromov, M, Positive curvature, macroscopic dimension, spectral gaps and higher signatures, in {\it Functional analysis on the eve of the 21st century}, Editors, S. Gindikin, et al, Birkh\"auser, Boston, 1996.



\bibitem[G-L]{gl2} Gromov, M, and Lawson, H., Positive scalar curvature and the Dirac operator on complete Riemannian manifolds, Publ. Math de IHES, {\bf 58} (1983) 83-196.


\bibitem[F]{F2}  Fraser, A., Fundamental groups of manifolds of positive isotropic curvature,  Ann. of Math. {\bf 158} (2003), 345--354.



\bibitem[M-W]{M-W} Micallef, M. and  Wang, M., Metrics with nonnegative isotropic
             curvature, {\em Duke Math. J.} {\bf 72} (1992), 649-672.


\bibitem[Ma]{ma} Massey, W., Algebraic Topology: An Introduction, GTM 56, Springer-Verlag, New York, 1984.


\bibitem[S-Y]{sy2} Schoen, R., and Yau, S. T., The existence of a black hole due to condensation of matter,  Comm. Math. Phys. {\bf 90} (1983) 575-579. 

\bibitem[Se]{s} Serre, J-P.,{\it Trees},  Springer-Verlag, Berlin, 1980.


\bibitem[St1]{s1} Stallings, J., {\it Group theory and three-dimensional manifolds}, Yale University Press, New Haven, 1971.


\bibitem[St2]{s} Stallings, J., On torsion-free groups with infinitely many ends, Ann. of Math. {\bf 88} (1968), 312-334.

\bibitem[W]{w} Wolfson, J., Manifolds of $k$-positive Ricci curvature, Proc. of the conference ``Variational problems in Riemannian geometry'', Leeds, UK, March 2009.



\end{thebibliography}
\end{document}